\documentclass{ifacconf}

\usepackage{natbib}

\usepackage{natbib}
\usepackage{amsmath,amssymb}
\usepackage{bm}
\usepackage{mathrsfs}
\usepackage{graphicx}
\usepackage{booktabs}
\usepackage{multirow}
\usepackage{siunitx}
\usepackage{algorithm}
\usepackage{algpseudocode}
\usepackage{tikz}
\usepackage{pgfplots}
\usepackage{xspace}
\usepackage{dsfont}
\usepackage{enumitem}
\usepackage{comment}
\usepackage{url}
\usepackage{dblfloatfix}
\usepackage{color}

\newtheorem{theorem}[thm]{Theorem}
\newtheorem{lemma}[thm]{Lemma}
\newtheorem{assumption}[thm]{Assumption}
\newtheorem{remark}[thm]{Remark}
\newtheorem{definition}[thm]{Definition}
\newtheorem{proposition}[thm]{Proposition}

\DeclareMathOperator*{\argmax}{arg\,max}
\DeclareMathOperator*{\col}{col}

\newcommand{\PP}{\mathbb{P}}
\newcommand{\EE}{\mathbb{E}}

\newcommand{\GP}[1]{\textcolor{blue}{#1}}

\begin{document}
\begin{frontmatter}

\title{Wasserstein Distributionally Robust Nash Equilibrium Seeking with Heterogeneous Data: A Lagrangian Approach\thanksref{footnoteinfo}}

\thanks[footnoteinfo]{%
The work of G. Pantazis and S. Grammatico has been partially supported by the ERC under project COSMOS (802348).
}

% Authors
\author[KTH]{Zifan Wang\thanksref{equal}}
\author[TUD]{Georgios Pantazis\thanksref{equal}}
\author[TUD]{Sergio Grammatico}
\author[ZAV]{Michael M.\ Zavlanos}
\author[KTH]{Karl H.\ Johansson}

% Equal-contribution footnote (if you still want the "$^*$" meaning)
\thanks[equal]{Equal contribution}

% Affiliations
\address[KTH]{Division of Decision and Control Systems, 
School of Electrical Engineering and Computer Science, 
KTH Royal Institute of Technology, 
SE-100 44 Stockholm, Sweden (e-mail: zifanw@kth.se, kallej@kth.se).}

\address[TUD]{Delft Center for Systems and Control, 
Delft University of Technology, 2628 CD Delft, 
The Netherlands (e-mail: G.Pantazis@tudelft.nl, S.Grammatico@tudelft.nl).}

\address[ZAV]{Department of Mechanical Engineering and Materials Science, Duke University, Durham, NC 27708, USA. (e-mail: michael.zavlanos@duke.edu}

\begin{abstract}                % Abstract of 50--100 words (IFAC guideline)
We study a class of distributionally robust games where agents are allowed to heterogeneously choose their risk aversion with respect to distributional shifts of the uncertainty. In our formulation, heterogeneous Wasserstein ball constraints on each distribution are enforced through a penalty function leveraging a Lagrangian formulation. We then formulate the distributionally robust game as a variational inequality problem, and show that under certain assumptions the original seemingly infinite-dimensional Nash equilibrium problem is equivalent to a multi-agent but finite-dimensional variational inequality problem with a strongly monotone mapping. Due to the inner maximization problem, it is however still challenging to calculate a distributionally robust Nash equilibrium. To this end, we design an approximate Nash equilibrium seeking algorithm and prove convergence of the average regret to a quantity that diminishes with the number of iterations, thus learning the desired equilibrium up to an a priori specified accuracy. Numerical simulations corroborate our theoretical findings. 
\end{abstract}

\begin{keyword}
Distributionally robust games \sep Wasserstein ambiguity sets \sep Nash equilibrium \sep Lagrangian methods \sep Variational inequalities
% (You can adjust to match the official IFAC keyword list.)
\end{keyword}

\end{frontmatter}

%===============================================================================

\section{Introduction}
\subsection{Motivation}
Modern engineering applications are increasingly operated by autonomous agents that interact strategically to optimize interdependent objectives subject to constraints, thus giving rise to a game-theoretic setting. Recent developments in control, optimization, and machine learning rely heavily on a data-driven multi-agent framework. As some examples,  power systems \cite{Sun2015, Molzahn2017} and networked robotics \cite{Bullo2009}, \cite{Cortes2017}, typically involve a high number of interacting decision-makers whose objectives are coupled through shared resources or performance constraints. At the same time, the environment these agents operate in is often uncertain, and its accurate probabilistic description is rarely known in practice. Thus, probabilistic models are often estimated from a finite collection of data, while statistical and algorithmic analyses are performed for the proposed methods. For uncertain multi-agent games \cite{fabiani_2022_automatica, Fele2019, Fele_CDC_2019,  Paccagnan_Campi_2019, PantazisFeleMargellos_Automatica2024} perform a detailed distribution-free probabilistic analysis for different game settings. From an algorithmic viewpoint, \cite{FranciGrammatico_TAC2021,FranciGrammatico_Automatica2022} design distributed schemes for stochastic generalized Nash equilibrium seeking. However, modeling uncertainty based on an empirical/ estimated distributions calculated on the basis of available samples and treated as ground truth can lead to lack of robustness towards even modest distribution shifts caused, for example, by changing operating conditions, sensor noise, seasonal demand variations or adversarial attacks. Thus, the decisions made are no longer (almost) optimal or as robust as expected.

Distributionally robust optimization (DRO) was developed precisely to tackle this mismatch between the true distribution and its estimated counterpart. In particular, instead of optimizing for one nominal distribution, DRO optimizes for the worst-case distribution within a considered ambiguity set. This equips the decision with quantifiable out-of-sample and finite-sample performance guarantees. Thanks to duality and optimal transport tools \cite{EsfahaniKuhn2018, AghassiBertsimas2006} even when the ambiguity set is a Wasserstein ball centered at the empirical distribution, the resulting reformulation can be  tractable, oftentimes formulated in the form of large linear programs under relatively mild assumptions. Classical DRO, however, is centralized, i.e.,  it assumes a single decision-maker and a single source of data. In contrast, in networked and learning-enabled systems, data are heterogeneous and private. Furthermore, as certain agents have different access to information/data, they might exhibit different levels of risk-aversion. A significant challenge that arises in several multi-agent systems is that agents can exhibit self-interested behaviour, wishing to minimize their own cost instead of a common social welfare cost. This framework gives rise to \emph{heterogeneous distributionally robust games}, where each agent hedges against their own distributional ambiguity while interacting strategically with others. 

In \cite{wang_icml22}, the authors develop a no-regret learning algorithm for multi-agent  convex games under CVaR-based risk measures.  The work in \cite{Wang_tac_2025}, studies convex games with heterogeneous feedback, while \cite{Wang_acc24} develops learning dynamics that converge to Nash equilibria in risk-averse games, where players optimize CVaR-type objectives. The work \cite{PantazisFranciGrammatico2023} relied on epigraphic reformulation methods to provide tractable reformulations for such games. When the agents are affected by a common ambiguity set based on shared data and ambiguity radius, the results in \cite{PantazisFranciGrammatico2023} recast the original DRNEP in a form which can provably be solved efficiently via centralized, decentralized or distributed algorithmic schemes. However, when ambiguity sets are private, this approach often leads to constraints which are private and coupling. An approach to fix this issue is to share the private constraints among agents. Unfortunately, this would introduce further cooperation among players, thus changing the original non-cooperative game with no clear connection between the equilibria of the two problems. To solve this problem, \cite{PantazisFranciGrammatico2023} propose a primal-dual algorithm inspired by the golden-ratio method, \cite{malitsky2020golden}, \cite{Franci2022_goldenratio}, to empirically converge to an equilibrium of the original problem. Following a different approach, \cite{Pantazis_dr_games_2024} focus on a particular class of DRNEPs and recast the problem into a variational inequality and propose efficient algorithmic methods to calculate an equilibrium. The recent paper in  \cite{Parsi2025} follows again a VI formulation, but instead proposes a gradient descent/ascent for which regret guarantees can be obtained under certain assumptions. \cite{lanzetti2025} introduce an alternative equilibrium concept based on optimal transport that protects against uncertainty in opponents' behavior, while \cite{fabiani_franci_2023} study distributionally robust generalized Nash games over  Wasserstein distributionally robust chance-constraints. 

\subsection{Main contributions}
Relative to the existing literature on Wasserstein distributionally robust Nash equilibrium problems and games our contributions are as follows:

\begin{enumerate}
    \item \textbf{Lagrangian formulation of DRNEPs with heterogeneous ambiguity sets:}
    To improve the computational efficiency in DRNE equilibrium seeking, while still maintaining certain distributional robustness properties, we propose a Lagrangian formulation of DRNEPs,
    where each agent enforces its Wasserstein ambiguity constraint via a penalty term in their cost, thus extending the formulation of the single-agent setting proposed in \cite{sinha2017certifying}. Our formulation enables heterogeneous and private ambiguity sets and
    risk-aversion levels without sharing constraints or data, while preserving the
    noncooperative nature of the game.
    We formalize the distributionally robust Nash equilibria (DRNE) for this novel setting and establish a connection with the DRNE of the original distributionally robust game. In our setting the agents' penalty parameters act as individual distributional robustness-tuning coefficients.

    \item \textbf{Monotone VI reformulation:}
    We reformulate the considered DRNEP as a finite-dimensional variational inequality problem and show that, under certain assumptions, we provide conditions ensuring (strong) monotonicity of the
    corresponding pseudo-gradient.
    We then design a DRNE seeking algorithm (Algorithm 1) in which agents
    draw data in an online fashion, and  prove regret guarantees of the averaged iterates
    towards an $\varepsilon$-DRNE, where $\varepsilon$ is set a priori. 
    We validate the proposed algorithm on a stochastic Nash-Cournot game with multiple risk-averse
    agents, and show efficient convergence to a neighbourhood of the DRNE with predicted accuracy and statistical properties of the algorithm with respect to out-of-sample performance. The proposed method can be adapted  for decentralized or distributed computation.
\end{enumerate}

\section{Problem Formulation}
\label{sec:1_problem}

We assume that each agent aims to minimize a stochastic cost function $f_i(x_i,x_{-i},\xi_i)$, where $x_i \in \mathbb{R}^{n}$ is the decision of agent $i$, $x_{-i} \in \mathbb{R}^{(N-1)n}$ denotes all other agents' decisions except for agent $i$, and $\xi_i \in \Xi_i$ is the individual uncertainty affecting agent $i$, where $\Xi_i$ is the corresponding support set. We assume that the agents do not have any knowledge of the support set $\Xi_i$ and/or the uncertainty of the probability distribution $\mathbb{P}_i$ that it follows. Instead, they have access to individual independent and identically distributed (i.i.d.) samples $\xi_i^{(k_i)}$, where $k_i=1,\dots, K_i$ and $K_i$ is the number of  samples for agent $i$. 
The empirical probability distribution for each agent $i$, denoted by $\hat{\mathbb{P}}_{K_i}$, is defined as
\begin{equation}
    \hat{\mathbb{P}}_{K_i} 
    \;:=\; \frac{1}{K_i} \sum_{k=1}^{K_i} \delta_{\xi_i^{(k)}},
\end{equation}
where $\delta_{\xi_i^{(k)}}$ denotes the Dirac measure at point $\xi_i^{(k)}$.

Finally, consider $X_i$ to be the local set in which the decision $x_i$ of each agent $i$ is confined. Consider a convex transport cost $c: \mathbb{R}^m \times \mathbb{R}^m \rightarrow \mathbb{R}$, which takes as input two realizations of random variables $z$ and $z'$ and returns a real number representing the distance between them. For example, consider $c(z,z') = \left\| z-z'\right\|_p^2$, where $\|\cdot\|_p$,  $p=1, 2, \dots, \infty $ denotes the $p$-Euclidean norm. Let $(\Xi_i,c)$ be a Polish metric space. For $p \ge 1$, the $p$-Wasserstein distance 
between two probability measures $\mathbb{Q}_1$ and $\mathbb{Q}_2$ on $\Xi_i$ 
is defined as
\begin{equation}
    W_p(\mathbb{Q}_1,\mathbb{Q}_2)
    \;:=\;
    \Bigg(
        \inf_{\pi \in \Pi(\mathbb{Q}_1,\mathbb{Q}_2)}
        \int_{\Xi_i \times \Xi_i} c(\xi,\zeta)^p \, \pi(\mathrm{d}\xi,\mathrm{d}\zeta)
    \Bigg)^{\!1/p},
\end{equation}
where $\Pi(\mathbb{Q}_1,\mathbb{Q}_2)$ denotes the set of all probability measures 
$\pi$ on $\Xi_i \times \Xi_i$ with marginals $\mathbb{Q}_1$ and $\mathbb{Q}_2$. Given a radius $\varepsilon_i > 0$, the Wasserstein ambiguity set for agent $i$ is given by:
\begin{align}
    \mathcal{P}_i
    \;:=\;
    \big\{
        \mathbb{Q} \in \mathcal{M}(\Xi_i)
        \;:\;
        W_p(\mathbb{Q},\hat{\mathbb{P}}_{K_i}) \le \varepsilon_i
    \big\},
\end{align}
where $\mathcal{M}(\Xi_i)$ denotes the set of all Borel probability measures on $\Xi_i$. A Wasserstein distributionally robust game is then defined as follows: 
\begin{equation}
\forall i \in \mathcal{N}: \min_{x_i \in X_i} \sup_{\mathbb{P}_i \in \mathcal{P}_i} 
\mathbb{E}_{\xi_i \sim \mathbb{P}_i}[J_i(x_i, x_{-i}, \xi_i)].
\tag{G}\label{DRGNEP}
\end{equation}
In other words, each agent solves their own optimization program that is coupled with other agents decisions through the cost function, while being risk-averse against the worst distribution within a Wasserstein ambiguity set $\mathcal{P}_i$.
To solve this problem we need to define a solution concept. The so-called \emph{distributionally robust Nash equilibrium} (DRNE) is provided in the following definition: 

\begin{definition} \label{def:DRNE}
A strategy profile $x^\star = (x_i^\star)_{i=1}^N \in \prod_{i=1}^N X_i$
is called a \emph{distributionally robust Nash equilibrium (DRNE)} of the game~$G$ if,
for every player $i \in \mathcal{N}$,
\begin{equation}
x_i^\star \in \arg\min_{x_i \in X_i} \;
\sup_{\mathbb{P}_i \in \mathcal{P}_i}
\mathbb{E}_{\xi_i \sim \mathbb{P}_i}
\!\left[f_i(x_i, x_{-i}^\star, \xi_i)\right].
\label{DRNE_def}
\end{equation}
\end{definition}
Note that the inner supremum in Definition \ref{def:DRNE} involves selecting an adversarial probability distribution that maximizes the cost of each agent. The seemingly infinite dimensional problem $(G)$ can be shown to admit tractable reformulations through duality and the use of epigraphic reformulations in \cite{PantazisFranciGrammatico2023} based on a methodology in the same spirit of \cite{EsfahaniKuhn2018}. In particular,  leveraging Kantorovich--Rubinstein duality, $(G)$ admits under mild assumptions the following reformulation:
\begin{equation}
\min_{x_i\in X_i,\;\lambda_i\ge 0}\;
\lambda_i \,\varepsilon_i
\;+\; \frac{1}{K_i}\sum_{k=1}^{K_i}
\sup_{\xi_i\in\Xi_i}\Big(f_i(x,\xi_i)-\lambda_i\|\xi_i-\xi_i^{(k)}\|_1\Big),
\label{eq:saddle}
\end{equation}
where \(\lambda_i\)  corresponds to a Lagrange multiplier of each Wasserstein ball constraint. It is then shown in \cite{PantazisFranciGrammatico2023} that if $(\tilde{x}, \tilde{\lambda})$ is an NE of this game, then $\tilde{x}$ is a DRNE OF $(G)$. However, computing a DRNE in this form
can be challenging as the number of agents and samples increases. To circumvent this challenge, in this work we propose an alternative reformulation based on the extended Lagrangian method. Specifically, instead of optimizing the Lagrange multiplier of the Wasserstein constraint, we instead impose it as a penalty in each cost. As such, each agent $i \in \mathcal{N}$ solves the following optimization problem:
\begin{align}\label{eq:obj0}
    \min_{x_i \in X_i} \max_{\PP_i}  \EE_{\xi_i \sim \PP_i} [f_i(x_i,x_{-i},\xi_i)] - \lambda_i W_p(\hat{\PP}_{K_i},\PP_i)
\end{align}
We call the collection of the interdependent optimization problems above game $G_L$, where the subscript $L$ denotes the Lagrangian relaxation of the original game $G$.
The following lemma then holds:
\begin{lemma} \label{lem:DR_reformulation}
The game $G_L$ in (\ref{eq:obj0}) can be equivalently reformulated as:
\begin{align}
  \forall i \in \mathcal{N}: \min\limits_{x_i \in X_i} \EE_{\hat{\xi}_i \sim \hat{\PP}_i} [\max\limits_{\xi_i\in \Xi_i}  f_i(x_i,x_{-i},\xi_i) - \lambda_i c(\xi_i,\hat{\xi}_i) ], \nonumber
\end{align} 
where $\EE_{\hat{\xi}_i \sim \hat{\PP}_i}[\cdot]$ denotes the expected value over the discrete empirical distribution.  \hfill $\square$
\end{lemma}
\emph{Proof}: The proof follows by adaptation of \cite{sinha2017certifying} for each agent $i \in \mathcal{N}$, given $x_{-i}$. \hfill $\blacksquare$

To further simplify the analysis, we consider the functions
\begin{align}
\hat{H}_i(x,\hat{\xi}_i)= \max_{\xi_i \in \Xi_i} h_i (x,\xi_i,\hat{\xi}_i) \nonumber 
\end{align}
and 
\begin{align}
h_i (x,\xi_i,\hat{\xi}_i)= f_i(x,\xi_i) - \lambda_i c(\xi_i,\hat{\xi}_i) \nonumber 
\end{align}
Thus, the optimization problem of each agent $i \in \mathcal{N}$ above takes the form:
\begin{align}\label{eq:obj1}
     \min_{x_i} \quad  & H_i(x) \nonumber \\
     &= \EE_{\hat{\xi}_i \sim \hat{\PP}_i } [\hat{H}_i(x, \hat{\xi}_i)] \nonumber \\
    &=\EE_{\hat{\xi}_i \sim \hat{\PP}_i } [\max_{\xi_i\in \Xi_i} \{ f_i(x,\xi_i) - \lambda_i c(\xi_i,\hat{\xi}_i^{j_i}) \}] \nonumber  \\
    &= \EE_{\hat{\xi}_i \sim \hat{\PP}_i } [\max_{\xi_i\in \Xi_i} h_i (x,\xi_i,\hat{\xi}_i) ],
\end{align}
Denote by $K$ the total number of drawn samples, i.e., $K=\sum_{i=1}^n K_i$. The distributionally robust NE for game $(G_L)$ is then defined as follows: 
\begin{definition}

%\begin{itemize}
%\item Consider game $G_K$. A point $\tilde{x}_K$ is called a strong distributionally robust Nash equilibrium (s-DRNE) if 
%\begin{align}
%\sup_{\mathbb{P}_i \in \mathcal{P}_i} &\mathbb{E}_{\xi_i \sim \mathcal{P}_i}[f_i(\tilde{x}_K, \xi_i)] \nonumber \\  
%&\leq \min_{x_i \in X_i} \sup_{\mathbb{P}_i \in \mathcal{P}_i} \mathbb{E}_{\xi_i \sim \mathcal{P}_i}[f_i(x_i, \tilde{x}_{-i, K}, \xi_i)]
%\end{align} for all $i \in \{1, \dots, n\}$
For the game $G_L$ in (\ref{eq:obj0}) with $\lambda=\col(\lambda_i)_{i=1}^N$, a point $x^\ast_\lambda$ is a distributionally robust Nash equilibrium (DRNE) if 
\begin{align} H_i(x_\lambda^\ast) 
&\leq \min_{x_i \in X_i} H_i(x_i, x^\ast_{-i, \lambda} ),
\end{align} for all $i \in \{1, \dots, n\}$. \hfill $\square$
%\end{itemize}
\end{definition}
\begin{remark}
If we select an appropriate penalty parameter $\lambda_i$ for each agent $i \in \mathcal{N}$ which coincides with the solution in the dual reformulation of $(G)$, then a DRNE of $(G_L)$ coincides with the standard DRNE. %However, computing these parameters through the dual formulation may lead to a nonmonotone mapping for heterogeneous ambiguity sets, whose equilibria can in general be hard to calculate. %This can hinder computational speeds and cause global optimality issues to emerge.
In this work we propose an algorithm that  efficiently computes an equilibrium of the alternative DRNEP in (\ref{eq:obj0}), while still be able to provably certify a principled level of distributional robustness for this point, by treating $\lambda_i$ as a penalty parameter. 
\end{remark}

\section{Distributionally Robust Nash Equilibrium seeking}

\subsection{Distributionally Robust Nash Equilbrium Problems as a Variational Inequality Problem}
In this section we recast the game $G_L$ in (\ref{eq:obj0}) as a variational inequality and establish the conditions that guarantee (strong) monotonicity of the corresponding mapping. We propose an equilibrium seeking algorithm to calculate efficiently a DRNE of (7). In particular, we propose Algorithm 1 which is based on an inexact proximal operation to calculate the inner supremum followed by a projected gradient descent for each agent $i \in \mathcal{N}$. We show that under certain conditions, convergence of the iterations to a neighbourhood of the equilibrium is guaranteed.

%\subsection{Strong Monotonicity Analysis}

We impose the following assumption: 
\begin{assumption} \label{assump:012}
For each agent $i \in \mathcal{N}$ it holds that:
\begin{enumerate}[label=\alph*)]
   \item   The decision set $X_i$ is closed and convex.
\item The cost function $f_i: \mathbb{R}^{nN} \rightarrow \mathbb{R}$ is differentiable and concave in $\xi_i$. \hfill $\square$
\end{enumerate}
\end{assumption} 
Based on Assumption \ref{assump:012}, the following result can then be obtained:
\begin{lemma} \label{lem:strong_concavity}
Under Assumption~\ref{assump:012}a), $h_i (x,\xi_i, \hat{\xi}_i)$ is $2\lambda_i$-strongly concave in $\xi_i$. \hfill $\square$
\end{lemma}
 Lemma \ref{lem:strong_concavity} makes the inner optimization problem easier to handle. Furthermore, let us impose the following assumption: 
\begin{assumption} \label{assump:2}
For any pair $x,y \in X=\prod_{i=1}^N X_i$, it holds that
\begin{align}
\left\| \nabla_{\xi_i} h_i (x,\xi_i,\hat{\xi}_i) -\nabla_{\xi_i} h_i(y,\xi_i,\hat{\xi}_i) \right\| \leq L_{\xi_i} \left\| x-y\right\|,
\end{align}
for all $\hat{\xi}_i$ and $i=1,\ldots,n$.
\end{assumption}
This assumption implies that the rate of change of the augmented Lagrangian $h_i$ with respect to the uncertain parameter $\xi_i$ is Lipschitz continuous with respect to $x$. To facilitate the analysis, we define
\begin{align}
\bar{M}_i(x,\hat{\xi}_i) = \argmax\limits_{\xi_i \in \Xi_i} [f_i(x,\xi_i) - \lambda_i c(\xi_i,\hat{\xi}_i)]. \nonumber 
\end{align}

%The above assumptions guarantee the well-behavedness of $H_i^{j_i}$, which is summarized in the following lemma.
The following Lemma then holds: 
\begin{lemma} \label{lemma:differentiable}
Given Assumptions ~\ref{assump:012} and \ref{assump:2}, $H_i$ is differentiable and satisfies
\begin{align}
    \nabla_{x_i} \hat{H}_i(x_i, x_{-i},\hat{\xi}_i) = \nabla_{x_i} f_i(x, \bar{M}_i(x,\hat{\xi}_i)),
\end{align}
and 
\begin{align}
    \left\| \bar{M}_i(x,\hat{\xi}_i) - \bar{M}_i(y,\hat{\xi}_i)\right\| \leq \frac{L_{\xi_i}}{2\lambda_i} \left\| x- y\right\|.
\end{align}
\end{lemma}
\emph{Proof}: The proof is an adaptation of Lemma 1 in \cite{sinha2017certifying} applied to the cost function of each agent $i \in \mathcal{N}$. \hfill $\blacksquare$ 

Consider now the pseudo-gradient mapping:
\begin{align}
F_K(x)=\text{col}((\nabla_{x_i}H_i(x))_{i=1}^N) \nonumber 
\end{align}
and $X=\prod_{i=1}^N X_i$. We define the  variational inequality problem \text{VI}$(F_K, X)$ as follows:
\begin{align}
\text{VI}(F_K, X): F_K(x^\ast)^\top(y-x^\ast) \geq 0 \text{ for all y }\in X. 
\end{align}

Then, the following proposition holds:
\begin{proposition}

    Consider Assumptions \ref{assump:012} and \ref{assump:2}. Then, the solution set of VI$(F_K, X)$ is equivalent to the set of DRNE of $G_L$ in (\ref{eq:obj0}). \hfill $\square$
\end{proposition}
\emph{Proof:} The proof can be obtained by applying Proposition 1.4.2 in \cite{FachineiPang} to  \text{VI}($F_K, X$). \hfill $\blacksquare$

The equivalence of the original DRNEP $G_L$ and \text{VI}($F_K, X$) allow us to leverage operator-theoretic methods to study and establish monotonicity properties of the original problem at hand. We impose the following assumption: 

\begin{assumption}\label{assump:3}
The gradient  $\nabla_{x_i} f_i(x,\xi_i)$ is $L_{x_i}$-Lipschitz continuous in $\xi_i$, i.e., 
\begin{align}
     \left\|\nabla_{x_i} f_i(x,\xi_i) - \nabla_{x_i} f_i(x,\xi_i')\right\| \leq L_{x_i} \left\| \xi_i - \xi_i'\right\|.
\end{align}
for all $i \in \mathcal{N}$.
\end{assumption}
Consider the pseudo-gradient 
\begin{equation}
G(x, \xi)=\col((\nabla_{x_i}f_i(x_i,x_{-i}, \xi_i)_{i=1}^N),
\end{equation}
where $\xi=\col((\xi_i)_{i=1}^N)$. We then consider the following assumption: 
\begin{assumption}\label{assump:4}
 $G$ is $\mu$-strongly monotone in its first argument for any $\xi_i\in \Xi_i$, $i \in \mathcal{N}$, i.e, for any decisions $x,y \in X$ and $\xi_i\in \Xi_i$ it holds that
\begin{align}
    (G(x, \xi)-G(y, \xi))^\top(x-y) \geq \mu \left\| x-y\right\|^2.
\end{align}
\end{assumption}

In the following result we establish the conditions on $\mu$ such that the mapping $F_K$ is also strongly monotone.

\begin{proposition}
If $\mu > \sqrt{n}\mu_{\xi}$,
where $\mu_{\xi}:= \max\limits_{i} \frac{L_{x_i} L_{\xi_i}}{2\lambda_i}$,
then $F_K$ is  $(\mu- \sqrt{n}\mu_{\xi})$-strongly monotone.
\end{proposition}

\emph{Proof:}
From Lemma~\ref{lemma:differentiable}, we have
\begin{align*}
    &\sum_{i=1}^{n} \EE_{\hat{\xi}_i \sim \hat{\PP}_i}  \Big\langle \nabla_{x_i} H_i(x,\hat{\xi}_i) - \nabla_{x_i} H_i(y,\hat{\xi}_i), x_i - y_i \Big\rangle \\
    & = \sum_{i=1}^{n} \EE_{\hat{\xi}_i \sim \hat{\PP}_i} \Big\langle \nabla_{x_i} f_i(x, \bar{M}_i(x,\hat{\xi}_i)  ) - \nabla_{x_i} f_i(y,\bar{M}_i(y,\hat{\xi}_i)), x_i - y_i \Big\rangle \\
    & = \sum_{i=1}^{n} \EE_{\hat{\xi}_i \sim \hat{\PP}_i} \Big\langle \nabla_{x_i} f_i(x,\bar{M}_i(x,\hat{\xi}_i)) - \nabla_{x_i} f_i(y,\bar{M}_i(x,\hat{\xi}_i)) \\
    &\quad + \nabla_{x_i} f_i(y,\bar{M}_i(x,\hat{\xi}_i)) - \nabla_{x_i} f_i(y,\bar{M}_i(y,\hat{\xi}_i)), x_i - y_i \Big\rangle \\
    &\geq \mu \left\| x-y\right\|^2 \\
    & + \sum_{i=1}^{n} \EE_{\hat{\xi}_i \sim \hat{\PP}_i} \Big\langle \nabla_{x_i} f_i(y,\bar{M}_i(x,\hat{\xi}_i)) - \nabla_{x_i} f_i(y,\bar{M}_i(y,\hat{\xi}_i)), x_i - y_i \Big\rangle \\
    &\geq \mu \left\| x-y\right\|^2 - \sum_{i=1}^{n} \frac{L_{x_i} L_{\xi_i} }{2\lambda_i} \left\| x-y\right\| \left\|x_i-y_i\right\| \\
    &\geq \mu \left\| x-y\right\|^2 -  \mu_{\xi} \left\| x-y\right\| \sum_{i=1}^{n}\left\|x_i-y_i\right\| \\
    &\geq \mu \left\| x-y\right\|^2 -  \sqrt{n}\mu_{\xi} \left\| x-y\right\|^2 \\
    &= (\mu- \sqrt{n} \mu_{\xi}) \left\|x-y\right\|^2,
\end{align*}
where the first inequality follows from Assumption~\ref{assump:4}, the second inequality follows since 
\begin{align*}
    &\left\| \nabla_{x_i} f_i(y,\bar{M}_i(x,\hat{\xi}_i)) - \nabla_{x_i} f_i(y,\bar{M}_i(y,\hat{\xi}_i)) \right\| \\
    &\leq L_{x_i} \left\|\bar{M}_i(x,\hat{\xi}_i) -\bar{M}_i(y,\hat{\xi}_i) \right\| \leq \frac{L_{x_i} L_{\xi_i}}{2\lambda_i} \left\| x - y\right\|.
\end{align*}
The statement follows if $\mu - \sqrt{n} \mu_{\xi}>0$. \hfill $\blacksquare$

\subsection{Average Convergence Guarantees}\label{sec:algorithm}

The following theorem provides average convergence guarantees for the proposed Algorithm 1:

\begin{algorithm}[t]
  \caption{Distributionally Robust NE seeking}\label{alg:DRG}
  \begin{algorithmic}[1]  % [1] = show line numbers
    \State \textbf{Input:} $\Xi_i$, parameter $\lambda_i$, and step size $\eta_i$, $i=1,\ldots,n$.
    \For{$t=1,2,\ldots$}
      \For{$i=1,\ldots,n$}
        \State Draw one sample $z_{i,t}$ from $\mathbb{P}_i$
        \State $\bar{z}_{i,t} \in \arg\max_{\xi_i \in X_i}
               \bigl(f_i(x_t,\xi_i) - \lambda_i c(\xi_i,z_{i,t}) + \epsilon_i\bigr)$
        \State $x_{i,t+1} = \operatorname{proj}_{X_i}\bigl(x_{i,t} - \eta_i \nabla_i f_i(x_t,\bar{z}_{i,t})\bigr)$
      \EndFor
    \EndFor
    \State \textbf{Output:} Sequence $x_t$
  \end{algorithmic}
\end{algorithm}

\begin{theorem} \label{thm:convergence}
Let Assumptions \ref{assump:012} and \ref{assump:2}, \ref{assump:3} hold. 
Then, the iterates $(x_t)_{t \in \mathbb{N}}$ of Algorithm 1 averagely converge to an $\epsilon$-DRNE, i.e.,  
\begin{align}
    \frac{1}{T} \sum_{t=1}^T \EE \left\| x_t 
 - x_{\lambda}^*\right\|^2 =\mathcal{O}(\frac{1}{\sqrt{T}} + \epsilon) 
\end{align}
where $\epsilon = \sum_{i=1}^n \epsilon_i$. \hfill $\square$

\end{theorem}

\emph{Proof:}
Assume that for any given $\bar{z}_{i,t}$ and $x_t$, 
\begin{align}
\left\| \nabla_{x_i} f_i (x_{i, t},x_{-i, t}, \bar{z}_{i,t}) \right\| \leq U, \nonumber 
\end{align}
for all $t \in \{1,\dots, T \}$. Furthermore, assume that the diameter of $\mathcal{X}$ is bounded by $D$. 
 Then, the following inequalities hold:
\begin{align}
    &\left\| x_{i,t+1} - x_{\lambda_i}^* \right\|^2  \nonumber \\
    &\leq \left\| x_{i,t} - \eta_i \nabla_{x_i} f_i (x_t,\bar{z}_{i,t}) - x_{\lambda_i}^* \right\|^2 \nonumber \\
    & = \left\| x_{i,t} -x_{\lambda_i}^* \right\|^2 - 2 \eta_i \langle x_{i,t} - x_{\lambda_i}^*, \nabla_{x_i} f_i (x_t,\bar{z}_{i,t}) \rangle \nonumber \\
    &\quad + \eta_i^2 \left\| \nabla_{x_i} f_i (x_t,\bar{z}_{i,t}) \right\|^2 \nonumber \\
    & =  \left\| x_{i,t} -x_{\lambda_i}^* \right\|^2 - 2 \eta_i \langle x_{i,t} - x_{\lambda_i}^*, \nabla_{x_i} \hat{H}_i(x_t,z_{i,t})\rangle + \eta_i^2U^2  \nonumber \\
    &\quad + 2 \eta_i \langle x_{i,t} - x_{\lambda_i}^*,\nabla_{x_i} f_i (x_t,\bar{z}_{i,t}) - \nabla_{x_i} f_i(x_t, \bar{M}_i(x,z_{i,t})) \rangle \nonumber \\
    &\leq  \left\| x_{i,t} -x_{\lambda_i}^* \right\|^2 - 2 \eta_i \langle x_{i,t} - x_{\lambda_i}^*, \nabla_{x_i} \hat{H}_i(x_t,z_{i,t})\rangle + \eta_i^2 U^2 \nonumber \\
    &\quad  + 2 \eta_i D L_{x_i} \left\| \bar{z}_{i,t} - \bar{M}_i(x,z_{i,t})\right\| \nonumber \\
    &\leq  \left\| x_{i,t} -x_{\lambda_i}^* \right\|^2 - 2 \eta_i \langle x_{i,t} - x_{\lambda_i}^*, \nabla_{x_i}  \hat{H}_i(x_t,z_{i,t})\rangle + \eta_i^2 U^2 \nonumber \\
    &\quad + 2 \eta_i D L_{x_i} \epsilon_i
\end{align}
The first inequality stems from the update step of the decisions in Algorithm 1 and the use of fundamental properties of the projection operator. The second equality is obtained by taking into account the property that 
\begin{equation}
\nabla_{x_i} \hat{H}_i(x_t,z_{i,t}) = \nabla_{x_i} f_i(x_t, \bar{M}_i(x,z_{i,t})). \nonumber
\end{equation}
The second inequality is obtained by applying the Lipschitz property of $\nabla_{x_i}f_i$, while the last one from the fact that the maximizer of the inner loop in Algorithm 1 calculates the solution from an initial point  $z_{i,t}$ with $\epsilon_i$ accuracy.
Taking expectations on both sides with respect to $z_{i,t}$, we have
\begin{align}
    &\EE \left\| x_{i,t+1} - x_{\lambda_i}^* \right\|^2  \nonumber \\
    &\leq \EE \left\| x_{i,t} -x_{\lambda_i}^* \right\|^2 - 2 \eta_i \langle x_{i,t} - x_{\lambda_i}^*, \nabla_{x_i}  H_i(x_t)\rangle + \eta_i^2 U^2 +  \nonumber \\
    &2 \eta_i D L_{x_i} \epsilon_i \nonumber \\
    &\leq \EE \left\| x_{i,t} -x_{\lambda_i}^* \right\|^2 - 2 \eta_i \langle x_{i,t} - x_{\lambda_i}^*, \nabla_{x_i}  H_i(x_t) - \nabla_{x_i}  H_i(x_{\lambda}^*) \rangle \nonumber \\
    &\quad  + \eta_i^2 U^2 + 2 \eta_i D L_{x_i} \epsilon_i \nonumber \\
\end{align}
where the last inequality follows from the Nash equilibrium first-order optimality condition.
Denote $\bar{\mu} = \mu - \sqrt{n} \mu_{\xi}$.
Summing up the above equation over $i=1,\ldots,N$ and setting $\eta_i=\eta$ yields
\begin{align}
    &\EE \left\| x_{t+1} - x_{\lambda}^* \right\|^2 \nonumber \\
    &\leq \EE \left\| x_{t} -x_{\lambda}^* \right\|^2 - 2 \eta \sum_{i=1}^N \langle x_{i,t} - x_{\lambda_i}^*, \nabla_{x_i}  H_i(x_t) - \nabla_{x_i}  H_i(x_{\lambda}^*) \rangle \nonumber \\
    &\quad  + N \eta^2 U^2 + 2 \eta D \sum_{i=1}^n L_{x_i} \epsilon_i \nonumber \\
    &\leq (1-2 \eta \bar{\mu} )\EE \left\| x_{t} -x_{\lambda}^* \right\|^2  + N \eta^2 U^2 + 2 \eta D \sum_{i=1}^n L_{x_i} \epsilon_i\nonumber \\
\end{align}
Rearranging the above equation and summing up over $t=0,\ldots,T-1$ yields
\begin{align}
    &\frac{1}{T} \sum_{t=1}^T \EE \left\| x_t - x_{\lambda}^*\right\|^2 \nonumber \\
    &\leq \frac{1}{2\eta \bar{\mu} T} \sum_{t=0}^{T-1} \big( \EE \left\| x_t - x_{\lambda}^*\right\|^2 - \EE \left\| x_{t+1} - x_{\lambda}^*\right\|^2\big) \nonumber \\
    &\quad + \frac{N \eta U^2}{2 \bar{\mu} } + \frac{D \sum_{i=1}^n L_{x_i} \epsilon_i} {\bar{\mu}} \nonumber \\
    &\leq \frac{D^2}{2\eta \bar{\mu} T} + \frac{N \eta U^2}{2 \bar{\mu} } + \frac{D \sum_{i=1}^n L_{x_i} \epsilon_i} {\bar{\mu}}.
\end{align}
Substituting $\eta = \frac{1}{\sqrt{T}}$ into the above equation completes the proof. \hfill $\blacksquare$

\section{Numerical Simulations}\label{sec:experiments}

\subsection{ Computational Efficiency of Algorithm 1}

In this section, we show the behaviour of Algorithm 1 on an illustrative example. 
Specifically, we consider risk-averse agents with the stochastic cost
\begin{align}
f_i(x,\xi_i)=%\bar{a}-
p(x_{-i})x_i+ c_i x_i+\xi_i x_i, \nonumber
\end{align}
where $x_i$ denotes the decision of agent $i$, $\xi_i\sim U(0,1)$ is a uniform random variable that represents uncertainty and $p(x_{-i})=a-\sum_jx_j$. Assume that each agent has $K_i$ samples of the random variable $\xi_i$, which we denote by $\hat{\PP}_i$. We set the accuracy of the inner loop $\epsilon = 10^{-3}$ in Algorithm 1.

\begin{figure}[h] 
\begin{center}
\centerline{\includegraphics[scale=0.75, trim=40pt 15pt 60pt 10pt, 
    clip]{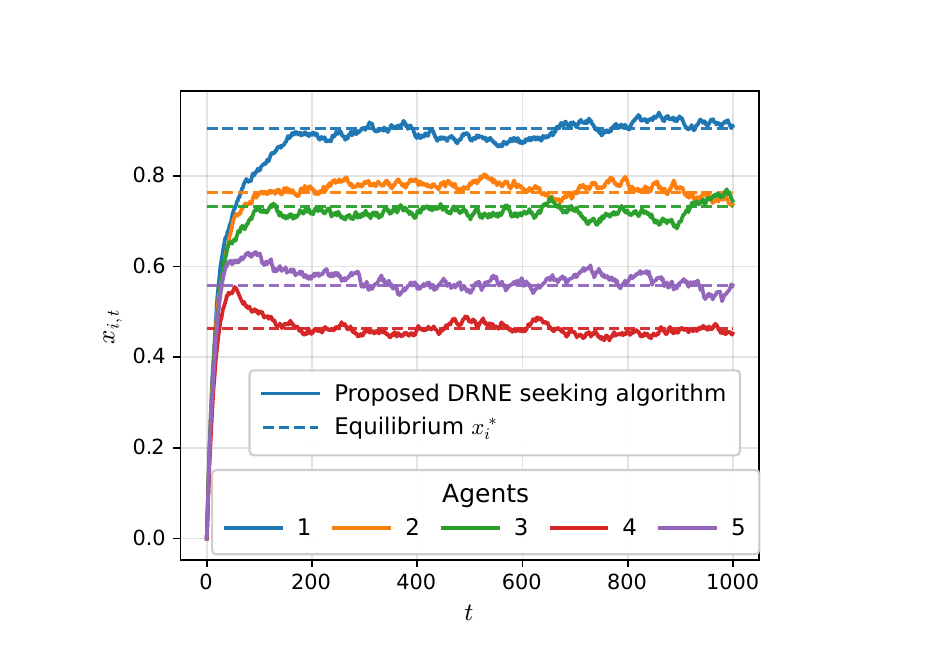}}
\caption{Decision iterates of each agent vs equilibrium solution. Note that the algorithm convergenges within a neighbourhood of the true DRNE values as stated in the provided guarantees.}
\label{fig_action}
\end{center}
\end{figure}

\begin{figure}[h] 
\begin{center}
\centerline{\includegraphics[scale=0.8]{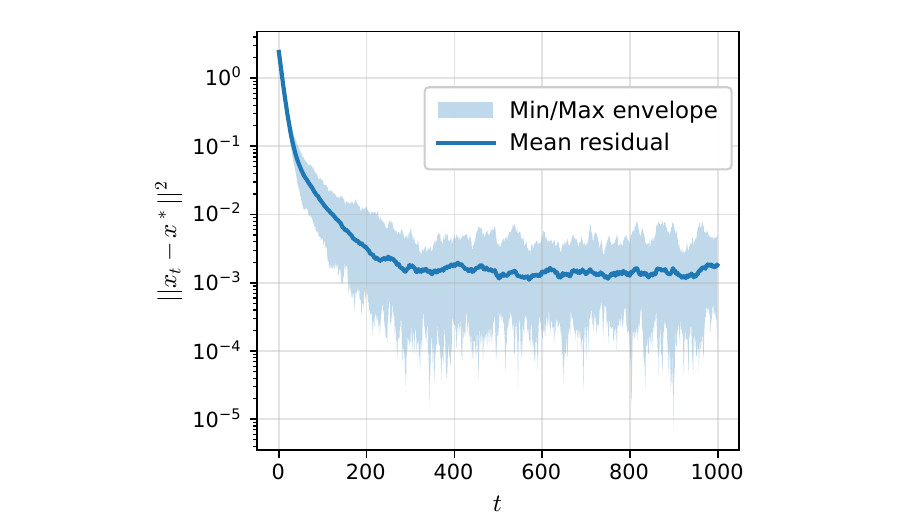}}
\caption{Convergence of the residual of Algorithm 1 using different multi-samples for training. The solid blue line corresponds to the mean among all residuals of the experiments, while the transparent blue cloud around the mean denotes the values that the residual can take across these repetitions.}
\label{fig_conv}
\end{center}
\vskip -0.2in
\end{figure}

As we can see from Figures 1 and 2, Algorithm 1 converges to the distributionally robust Nash equilibrium with the predefined accuracy. Furthermore, the proposed algorithmic scheme is fast and can be used for decentralized or distributed design. Figure \ref{fig_action} illustrates the agents' iterates for $N=5$. Note that the iterates of each agent converge to an $\epsilon$-neighbourhood of the equilibrium solution, thus validating numerically the result of Theorem \ref{thm:convergence}. Figure \ref{fig_conv}  shows the residual between the equilibrium and the iterates for Algorithm 1 for experiments where different samples are used for learning the equilibrium. The solid blue line represents the mean among all the experiments, while the transparent clouds represent the empirical variance across the experiments.
 \subsection{Statistical performance}
Now we wish to test the out-of-sample performance of Algorithm 1 based on the risk-aversion of the individual agents. Note that smaller values of $\lambda_i$ in Algorithm 1 imply that agents are not penalized significantly to follow the empirical distribution. %Thus, they are more risk-aware towards deviations from the estimated distribution.
Larger values of $\lambda_i$ imply that agent $i$ puts more emphasis on following and trusting the empirical distribution and is not allowed to deviate much from it. As such, the inverse term $\rho_i= \frac{1}{\lambda_i}$ can be viewed as a parameter that shares similaries with the Wasserstein radius in standard DRNEPs as an increase in its value implies a higher risk aversion towards probabilistic shifts. 

\begin{table*}[t]
\centering
\begin{tabular}{c l l l}
\hline
\textbf{Case} &
\multicolumn{1}{c}{Scenario 1} &
\multicolumn{1}{c}{Scenario 2 } &
\multicolumn{1}{c}{Scenario 3 } \\
\hline
1 &
$\mathbf{\rho}=[0.05,\ 0.075,\ 0.10,\ 0.125,\ 0.15]$ &
$\mathbf{\rho}=[2.0,\ 0.075,\ 0.10,\ 0.125,\ 0.15]$ &
$\mathbf{\rho}=[2.0,\ 2.0,\ 0.10,\ 0.125,\ 0.15]$ \\
2 &
$\mathbf{\rho}=[0.20,\ 0.30,\ 0.40,\ 0.50,\ 0.60]$ &
$\mathbf{\rho}=[4.0,\ 0.30,\ 0.40,\ 0.50,\ 0.60]$ &
$\mathbf{\rho}=[4.0,\ 4.0,\ 0.40,\ 0.50,\ 0.60]$ \\
3 &
$\mathbf{\rho}=[1.20,\ 1.30,\ 1.40,\ 1.50,\ 1.60]$ &
$\mathbf{\rho}=[8.0,\ 1.30,\ 1.40,\ 1.50,\ 1.60]$ &
$\mathbf{\rho}=[8.0,\ 8.0,\ 1.40,\ 1.50,\ 1.60]$ \\
\hline
\end{tabular}
\end{table*}

\begin{figure}[t]
    \centering
    \includegraphics[width=\columnwidth]{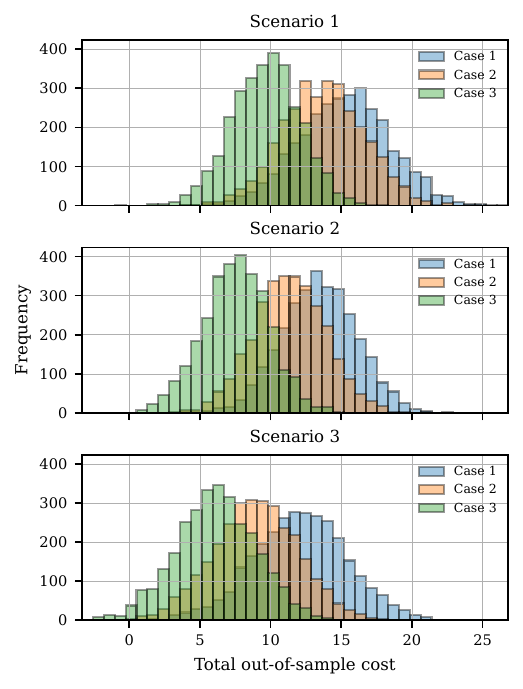}
    \caption{Histogram of out-of-sample total cost for the following scenarios 1) No agent shows signficantly higher conservatism compared to the others; 2) Agent 1 shows higher levels of conservatism 3) Agents 1 and 2 show higher levels of conservatism compared to other agents.}
    \label{fig:oos-hists}
\end{figure}

To further test the out-of-sample performance of    Algorithm 1 we consider that $\xi_i$ follows a Gaussian distribution of mean $\mu_i$ and variance $\sigma^2_i$. We denote the collections of means and standard deviations for the individual distributions affecting the agents by $\mathbf{\mu}=\{\mu_i\}_{i=1}^5$ and $\mathbf{\sigma}=\{\sigma_i\}_{i=1}^5$. We set these values at $\mu=(0, 0.3, 0.6, 0.9, 1.2)$ and $\sigma=(1, 1.2, 1.5, 1.8, 2)$. To test the out-of-samples performance we consider that the mean and standard deviation vectors are perturbed by $\delta \mu=(0.5, -0.4, 0.6, -0.5, 0.7)$ and $\delta \sigma=(0.8, -0.6, 0.9, -0.7, 1)$ and draw $K=(K_i)_{i=1}^5=(20, 15, 10, 8, 6)$ and a total number of $K_{test}=3000$ random samples vectors for testing. Based on the drawn samples we calculate the total cost of the population per sample vector realization. Figure \ref{fig:oos-hists} illustrates the out-of-sample performance for three different scenarios of parameter choices $\rho_i$, $i \in \mathcal{N}$ as shown in Table 1. Scenario 1 involves all agents selecting different risk-aversion parameters $\rho_i$ but on a similar scale. In Scenario 2, agent 1 shows significant risk-aversion compared to the other agents, while in Scenario 3 both agents 1 and 2 are significantly more risk-averse compared to the other three agents. Based on Figure \ref{fig:oos-hists}, Scenario 3 shows better out-of-sample perfomance for the population cost overall, as evidence by the shift of the histogram towards the left. In terms of out-of-sample performance, Scenario 2 is still overall better than Scenario 1. This implies that due to the interconnection of the agents through the coupling in the cost functions, the risk-aversion of just some players will affect the out-of-sample performance of the entire population. For cases 1, 2, 3 (see Table 1) the figure shows that if all agents increase their risk-aversion, distributionally more robust equilibria are obtained.

\section{Conclusion}
\label{sec:conclusion}

We have introduced a Lagrangian formulation for data-driven Wasserstein distributionally robust games that lets agents independently tune their risk aversion while preserving their strategic self-interested behavior. By reformulating the resulting game as a finite-dimensional variational inequality, we derived explicit conditions for strong monotonicity of the pseudo-gradient and proposed a projected gradient--type learning algorithm whose averaged iterates converge to an approximate distributionally robust Nash equilibrium. Numerical experiments on a game with heterogeneous risk-averse players confirm the established convergence properties and highlight how the Lagrangian penalties affect distributional robustness for the population.  The next step is to tackle games with shared constraints paving the way towards an algorithmic theory for distributionally robust generalized Nash equilibria.

\bibliography{ifacconf} 

% or whatever the template says
%\bibliography{references} 

\end{document}